\documentclass[11pt]{amsart}
\usepackage{amsmath,amsthm, amscd, amssymb, amsfonts}
\usepackage{epsfig}
\input xy
\xyoption{all}
\usepackage[dvips, dvipsnames, usenames]{color}

\newcommand{\Ku}{{\mathbf K}}

\newcommand{\Z}{{\mathbb Z}}
\newcommand{\Na}{{\mathbb N}}

\newcommand{\R}{{\mathbb R}}

\newcommand{\id}{\mbox{\rm id\,}}

\newcommand\Hom{\operatorname{Hom}}

\theoremstyle{plain}

\theoremstyle{definition}

\theoremstyle{remark}

\def\pf{\begin{proof}}
\def\epf{\end{proof}}

\theoremstyle{remark}

\begin{document}
\title[the beginnings of the theory of Hopf algebras]
{The beginnings of the theory of Hopf algebras}
\author[andruskiewitsch]{Nicol\'{a}s Andruskiewitsch}
\address{Facultad de Matem\'{a}tica, Astronom\'{i}a y
F\'{i}sica, Universidad Nacional de C\'{o}rdoba, CIEM - CONICET
\\ \newline\indent (5000) Ciudad Universitaria \\ C\'{o}rdoba
\\Argentina} \email{andrus@famaf.unc.edu.ar}
\thanks{The first author
was  supported by Mincyt (Cba), CONICET, Foncyt and
Secyt (UNC). He thanks Dominique Flament for inviting him to deliver a course at the
\emph{\'Ecole d'\'et\'e d'histoire des math\'ematiques} (Brasilia, UnB, mars 2008). Preliminary material for
this paper was collected for that occasion. }
\author[ferrer]{Walter Ferrer Santos}
\address{Facultad de Ciencias\\\newline\indent Universidad de la
  Rep\'ublica \\\newline\indent Igu\'a 4225, Montevideo, Uruguay}\email{wrferrer@cmat.edu.uy}
\thanks{The second author was partially supported by Pedeciba (Ur) and Foncyt (Ar)}
\thanks{This version contains corrections to the one published in Acta
  Applicandae Mathematicae, some of them
suggested by a referee in a report that the authors received
after publication. }
\begin{abstract} We consider issues related to the
  origins, sources  and initial motivations of the theory of Hopf algebras.
  We consider the two main
  sources of primeval development: algebraic topology and algebraic
  group theory. Hopf algebras are named from the work of Heinz Hopf in
  the 1940's. In this note we trace the infancy of the subject back to
  papers from the 40's, 50's and 60's in the two areas mentioned
  above. Many times we just
describe -- and/or transcribe parts of --
  some of the relevant original papers on the subject.
\end{abstract}

\maketitle

\section{Introduction}
In this note, we address the following questions: Who introduced
the notion of Hopf algebra? Why does it bear this name?

Succinctly,
the answers are:
\renewcommand{\descriptionlabel}[1]
{\hspace{\labelsep}\textsf{#1}}
\begin{description}
    \item[$\ast$] The first formal definition of Hopf algebra -- under the name
      of \emph{hyperalgebra} --
    was given by
    Pierre Cartier in 1956, inspired by the work of Jean
    Dieudonn\'e on algebraic groups in positive characteristic.

    \item[$\ast$] The expression \emph{Hopf algebra}, or more precisely,
      \emph{alg\`ebre de
Hopf}, was coined by Armand Borel in 1953, honoring the
foundational work of Heinz Hopf.
\end{description}

That is, we have two different strands, one coming from the theory
of algebraic groups and the other from algebraic topology, both
interwoven in the early 60's. But the story is not so
simple: the notion introduced by Cartier is not exactly the same
as what we use today and it took a time to arrive at the present
definition. Also, what did Borel mean when speaking about Hopf
algebras? We discuss these topics quoting from the original
sources.


\medbreak This paper has three sections besides this introduction.
The second is centered
around those aspects of algebraic topology that played an
important role in the development of Hopf algebra theory. We
organized our expositions in terms of the different authors who
in our opinion were the main contributors. In particular we
discuss the work of Heinz Hopf, Armand Borel,
Edward Halpern, John Milnor and John C. Moore.

In the third section called algebraic groups we proceed
similarly and center our attention on the input of algebraic
group theory in Hopf theory. We discuss in particular the relevant
work of the following mathematicians: Jean Dieudonn\'e, Pierre
Cartier, Gerhard Hochschild, Dan Mostow, Bertram Kostant and
Georg Kac. In the fourth section we make a few general remarks
concerning later developments. \medbreak

We finish this introduction by saying some words about notations
and prerequisites. We assume that the reader has --besides a working
knowledge in Hopf algebra theory-- a certain degree of
familiarity with the basic vocabulary on: algebraic topology,
algebraic groups, Lie theory and (co)homology of Lie algebras,
etc. In Hopf algebra theory we use freely the standard modern
terminology as it appears for example in Moss Sweedler's book
\cite{Sw}.

\medskip
One of the goals of the present article is to
trace back the formal definition of a Hopf algebra and to describe its
evolution in time. The definition we know and use today is included in
the subsection on the work of Bertram Kostant, see page
\pageref{kostant}. Compare also with \cite[vol. 1, p. 511]{dg1}.

\medskip
{\sc Acknowledgments.} We thank Marcelo Aguiar, Gerhard
Hochschild, Jim Humphreys, Ber\-tram Kostant, Susan Montgomery and
Moss Sweedler, for commenting on an initial version of this paper.
However, any inaccuracy or omission is the responsibility of the
authors, who are not professional historians and had to send the
printers a first version of the paper under the pressure of a very
tight deadline.

\section{Algebraic topology}
\subsection{Heinz Hopf} In \cite{hopf} the author considers what he
calls a $\Gamma$--manifold -- later called a Hopf manifold by
Borel in \cite{borel53} and currently simply an $H$--space. They
are spaces equipped with a product operation i.e. a continuous
function $M \times M \rightarrow M$. Clearly this function induces
by functoriality a homomorphism from the cohomology ring of $M$
 -- called $H$ -- into the cohomology of $M \times M$, i.e. into $H
\otimes H$.  In the above mentioned paper Hopf shows that the
existence of this additional map from $H$ into $H \otimes H$ --
that is compatible with the cup product --  imposes strong
restrictions on the structure of $H$, and from these restrictions
he deduces some important topological results.  In accordance with
W. Hurewicz's review of \cite{hopf} (appearing in MathSciNet:
MR0004784, 56.0X):

\medskip
 {\em This [...]\/\,is not only a far reaching generalization of the well--known
  results of Pontrjagin about Lie groups, but in addition it throws an
  entirely new light on these results by emphasizing their elementary
  topological nature\ldots}
\medskip

The methods introduced by Hopf in this article, were soon applied
to other cohomology theories -- see for example \cite{leray} and
\cite{lefschetz}.

In a similar direction Hans Samelson, a close collaborator of
Hopf, following the ideas of the latter but working with the
homology instead of the cohomology (see \cite{sam}), considers two
products (the ``intersection'' and the ``Pontrjagin'' product),
as an alternative to the consideration {\em \`a la} Hopf of a
product and a coproduct on the cohomology.

 It is worth mentioning another more algebraic derivation of the
methods introduced by Hopf. Using de Rham theory, the study of the
real (co)homology of a Lie group -- a very special case of a Hopf
manifold (or $H$--space)-- can be reduced to the study of the {\em
algebraic} (co)homology of the Lie algebra ${\mathfrak g}$ of $G$.
J.L. Koszul in \cite{koszul} considers --for a {\em general} Lie
algebra $\mathfrak g$-- many problems related to the cohomology
and homology of Lie algebras and in particular deals with the
version for Lie algebras of the results of \cite{hopf}.

This generalization is achieved in the form of an identification
--under the hypothesis of reductivity-- of the cohomology as well
as the homology of $\mathfrak g$, to the exterior algebras of its
respective primitive elements. In the proof of this result, a
crucial technical role is played by the map $X \mapsto (X,X):
\mathfrak g \rightarrow \mathfrak g \times \mathfrak g$ that is
(c.f. C. Chevalley MathSciNet MR0036511) {\em \ldots an
  algebraic counterpart to the mapping $(s,t) \mapsto st$ of $G \times
  G$ into $G$ \ldots}

\subsection{Armand Borel} In the important paper by A. Borel
\cite{borel53}, the author studies the homology of a principal
fiber bundle and the application to the homology of homogeneous
spaces. Concerning this problem he mentions many mathematicians
who studied it before: H. Cartan, Chevalley, Koszul, Weil, around
1950, with predecessors in the work of E. Cartan, Ehresmann,
Pontrjagin (see the references in \cite{borel53}), Hopf and
Samelson (see \cite{hopf}, \cite{hopfsam} and \cite{sam}).

In the introduction to the paper,
when mentioning the results he needs from the work of
Hopf he uses the expression `Hopf algebra', or more precisely `alg\`ebre de
Hopf', that appears first in \cite[p.116]{borel53}:

\medbreak\noindent
\emph{\ldots\/la structure d'une alg\`ebre de Hopf (c'est
\`a dire v\'erifiant les conditions de Hopf)\ldots}
\medbreak

As we mentioned, Borel means the conditions in the paper
\cite{hopf} that are satisfied by the cohomology ring of an
$H$--space and that concern the existence of an additional
operation of comultiplication on the cohomology ring.

In the introduction the author mentions Hopf algebras over
fields of positive characteristic and in the same page remarks about
a ``structure
theorem'':

\medbreak\noindent
{\em On peut exprimer le r\'esultat
  obtenu en disant qu'une alg\`ebre de Hopf est toujours isomorphe
  \`a un produit tensoriel gauche d'alg\`ebres de Hopf \`a un
  g\'en\'erateur}.
\medbreak

Later, in the second chapter of the paper called ``Le th\'eor\`eme
de Hopf'' -- more precisely in \cite[p. 137]{borel53} -- he
considers an algebra $H$ with the following properties: it is
graded by $\mathbb N_0$, with finite-dimensional homogeneous
components; and it is anti-commutative (or super-commutative in
present terminology). He denotes by $D$ the degree. He formally
defines:

\medbreak \noindent  \emph{$H$ est une alg\`ebre de Hopf s'il
existe un homomorphisme d'alg\`ebres gradu\'ees $f:H \to H\otimes
H$ et deux automorphismes $\rho$ et $\sigma$ de $H$ tels que pour
tout $h\in H$ homog\`ene l'on ait:
\begin{equation}\label{borel6.1}\tag{6.1}
\begin{aligned}
f(h) = & \rho(h)\otimes 1 + 1\otimes \sigma(h) +x_1\otimes y_1 +
\dots + x_s\otimes y_s\\ &(0 < Dx_i < Dh).\end{aligned}
\end{equation}}
\medbreak

The author then shows that without loss
of generality it can be assumed that $\rho = \sigma = \id$.

Notice that the comultiplication map $f$ of a `Hopf algebra'
in the sense of Borel is neither coassociative
nor admits an antipode; and although it is not
explicitly said, it does admit a counit. This comes from the fact that the
original geometric multiplication considered by Hopf on the base
space, was not assumed to be
associative or to have an inverse, but in general it was assumed to possess a
``homotopy identity''.

Then the author proves the ``structure theorem'' mentioned before
-- Theorem 6.1 -- that is stated informally on p. 138 as follows:

\medbreak
\noindent{\em Ce th\'eor\`eme affirme en somme que $H$ est
  l'alg\`ebre associative engendr\'ee par les $x_i$, [..]
il montre aussi qu'une alg\`ebre de Hopf peut
  toujours \^etre envisag\'ee comme produit tensoriel (gauche)
  d'alg\`ebres de Hopf \`a un g\'en\'erateur}.
\medbreak

It is worth noticing that this structure theorem, that generalizes
the corresponding result due to Hopf, gives information only about
the multiplicative structure of the ``Hopf algebra''. Information
related to the coproduct is obtained in \cite{sam}, and also in
work of Leray \cite{leray2}. As we shall see in Subsection
\ref{subsection:Milnor--Moore}, results of this type were later
generalized and unified in the work of Milnor and Moore
\cite{milnormoore}.

The author draws immediately some important consequences of the
structure theorem above, e.g. in Chapter III p. 147 he proves
Proposition 10.2 -- in Borel's notation $K_p$ stands for the
finite field with $p$ elements, ${\bf{O}}(n)$  for the orthogonal
group of the $n$--space and ${\bf S}_n$ for the $n$--sphere:

\medskip
\noindent{\em \ldots\/ ${\bf{V}}_{n,n-q}= {\bf{O}}(n)/{\bf{O}}(q)$ ($1 \leq q \leq
  n$)\ldots  Soit $\bar{n}$, (resp. $\bar{q}$), le plus grand, (resp. le
  plus petit), entier impair $\leq n$ (resp. $\geq
  q$). Alors ${\bf V}_{n,n-q}$ a pour les coefficients $K_p$ ($p \neq 2$)
  m\^eme alg\`ebre de cohomologie que le produit:
\[ {\bf{S}}_{2\bar{n}-3} \times {\bf{S}}_{2\bar{n}-7} \times \cdots \times {\bf{S}}_{2\bar{q}+1}\]
 multipli\'e encore par ${\bf S}_{n-1}$ si $n$ est pair, par ${\bf
   S}_q$ si $q$ est pair}.
\medskip

It seems remarkable that in the seminal paper by A. Borel: {\it
Groupes lin\'eaires alg\'ebriques} \cite{borel56}, that was one of
the foundational stones of the modern theory of algebraic linear
groups for which the algebra of polynomial functions is the
archetype of a commutative Hopf algebra, neither the term `Hopf
algebra', nor any consideration about the Hopf structure of its
algebras of coordinate functions appear explicitly.

\subsection{Edward Halpern} In \cite{halpern1} the author presents
some results that appeared in his thesis under the direction of
S. MacLane that develops further the theory started in \cite{hopf} and
\cite{borel53}. See also \cite{halpern2} and \cite{halpern3} of the same
author.

In the introduction to \cite{halpern1} p. 2 the author says:

\medskip
\noindent {\em Thus, in both homology and cohomology for
$H$--spaces one is led to consider a triple consisting of a module
$H$, a product $H \otimes H \rightarrow H$ and a ``coproduct'' $H
\rightarrow H \otimes H$, with product and coproduct compatible.
Such a structure is called a \underline {hyperalgebra} -- the term
hyperalgebra was suggested by S. MacLane [and had previously been
used by Dieudonn\'e]\,\footnote{Citation from \cite{halpern2}}.
Both Hopf and Pontrjagin algebras are special cases of
hyperalgebras}.

\medskip

With the names of Hopf and Pontrjagin algebras the author is
referring to the following situations. If one starts with an
$H$--space $X$ (endowed with a continuous map $\Delta: X \times X
\rightarrow X$) and a homotopy unit, one can consider the induced
maps $\Delta_{*}: H_* \otimes H_{*} \rightarrow H_{*}$ and
$\Delta^{*}: H^* \rightarrow H^* \otimes H^{*}$  on homology and
cohomology respectively. Following the nomenclature of
\cite{borel53}, if $\cup$ denotes the cup--product, the triple
$(H^*, \cup, \Delta^*)$ was being called by some authors a {\em
Hopf algebra}. Dually the triple $(H_{*}, \Delta_{*}, \cup_{*})$
where $\cup_{*}$ is the dual of the cup--product was called a
Pontrjagin algebra in \cite{borel54}. It is a standard fact that
the maps cup and $\Delta^*$  as well as the maps $\cup_{*}$ and
$\Delta_{*}$ are compatible.

Hence, in the nomenclature of the author the structure maps of a
hyperalgebra need not be associative or coassociative,
but in general he supposes the underlying spaces to be
graded and the operations to preserve the grading, and to be unital
and counital. In the case of
what he calls a Hopf algebra the product is associative and in the case
of the so called Pontrjagin algebra the coproduct is coassociative.

Later in Chapter 1, p. 4, the hyperalgebra is said to be
associative if the product is associative and the coproduct is
coassociative.

In the later paper \cite{halpern2} the author continues to apply the name of
``hyperalgebra'' to the objects considered above, but in
\cite{halpern3} he switches to the name of ``Hopf algebra'' that was
becoming more and more popular after the work of Milnor and Moore.

\subsection{John W. Milnor and John C. Moore}
\label{subsection:Milnor--Moore} The  paper of  John W. Milnor and
John C. Moore that we consider appeared in print in March 1965,
but in accordance with the review MathSciNet: MR0174052 written by
I. M. James:

\medbreak \noindent {\em $\cdots$ [A] version of this
  paper circulated some ten years ago. There have been many
  improvements, although the published version is basically the
  same\,\footnote{The author of the review was probably referring
to J. Milnor and J. C. Moore, {\em On the structure of Hopf
algebras},\/ Princeton Univ., Princeton, N.J., 1959 (mimeographed)
or some predecessor.}}. \medbreak

The article has eight sections, in Section 4 the authors
introduce the concept of Hopf algebra, that coincides with the usual
definition of {\it graded bialgebra} that we use today. In the case
that either the product or the coproduct are not associative, the
authors speak of a {\em quasi Hopf algebra}.

\medbreak The authors work in the category of $\mathbb N_0$-graded
modules over a fixed commutative ring $K$. The definition of Hopf
algebra appears as Definition 4.1 on page 226:

\medbreak \noindent{\sc Definition}.\,{\it A Hopf algebra over $K$
is a graded
  $K$--module together with morphisms of graded $K$--modules
\begin{align*} \varphi: A \otimes A \rightarrow A,& \hspace{.3cm} \eta: K
    \rightarrow A \\ \Delta: A \rightarrow A \otimes A,& \hspace{.3cm}
    \varepsilon: A \rightarrow K
\end{align*}
such that
\begin{enumerate}
\item \,$(A,\varphi,\eta)$ is an algebra over $K$ with augmentation
  $\varepsilon$,
\item \,$(A,\Delta,\varepsilon)$ is a coalgebra over $K$ with augmentation
  $\eta$, and
\item the diagram
\[\xymatrix{A \otimes A \ar[r]^{\varphi}\ar[d]^{\Delta \otimes \Delta}
& A \ar[r]^{\Delta} & A \otimes A \\
A \otimes A \otimes A \otimes A \ar[rr]^{A \otimes T \otimes A} && A
\otimes A \otimes A \otimes A \ar[u]_{\varphi \otimes \varphi}}
\]
is commutative.
\end{enumerate}}
\medbreak

In the definition above -- and throughout this paper -- the authors
use the notation $A$ to indicate at the same time an object as
well as the identity map on that object, moreover the map $T$
stands for the usual braiding in the category of $\Na_0$--graded
modules.

All in all, what the authors call a Hopf algebra, is what today
would be called a bialgebra in the symmetric category of graded
modules over a commutative ring $K$. But connected Hopf
algebras in their sense, that is those with one-dimensional
0-degree space, do have an antipode as explained below on page
\pageref{antipode}.

Later, in Section 5, the authors give the structure Theorem 5.18,
today known as the Milnor--Moore theorem:

\medbreak \noindent{\it Over a field of characteristic zero, the
category of graded Lie algebras is isomorphic with the category of
primitively generated Hopf algebras.} \medbreak

In Section 7, {\em Some classical theorems}, a proof of a general
theorem on the structure of graded--commutative connected Hopf
algebras is presented. This theorem generalizes the array of
structural theorems mentioned before due to Hopf, Samelson, Leray,
Borel, etc. Along the way the authors prove the following
classical results: Theorem 7.5 (Leray), Theorem 7.11 (Borel),
Theorem 7.20 (Samelson, Leray). The main result proved in this
section reads as follows (Theorem 7.16, page 257).

\medbreak \noindent {\em If $A$ is a connected primitively
generated Hopf algebra over the
  perfect field $K$, the multiplication in $A$ is commutative, and the
  underlying vector space of $A$ is of finite type, then there is an
  isomorphism of Hopf algebras of $A$ with $\bigotimes_{i \in I}A_i$
  where each $A_i$ is a Hopf algebra with a single generator $x_i$}.
\medbreak

In the case that the hypothesis of primitive generation is not valid,
all that can be obtained is Borel's theorem stating the isomorphism
of $A$ with $\bigotimes_{i \in I}A_i$ as algebras.

\medbreak In the Appendix and as an illustration of the powerful
methods developped, using a combination of results from Hopf
theory and algebraic topology the authors  prove the following
very precise result for the Hopf algebra structure of the homology
ring of a Lie group -- formulated more generally for certain
$H$--spaces:

\medbreak \noindent{\em If $G$ is a pathwise connected homotopy
associative
  $H$--space with unit, and $\lambda: \pi(G;K) \rightarrow H_{*}(G;K)$
is the Hurewicz morphism of Lie algebras, then the induced
morphism $\widetilde{\lambda}:U(\pi(G;K)) \rightarrow H_{*}(G;K)$
is an isomorphism of Hopf algebras}. \medbreak

This is proved using the fact that for an associative connected
$H$--space the homology is commutative and connected, thus -- as
was proved before -- it is also primitively generated. As the
image of the map $\lambda$ is exactly the set of primitive
elements of $H_{*}(G;K)$ the result follows.

\medbreak

The article we are considering had a wide and deep influence in
the development of the subject, and maybe should  be thought of at
the same time as a culmination of the topological line of work
initiated by Hopf et.~al., as well as the launching platform of an
independent new area in the realm of abstract algebra. Moreover it
was one of the basic references on the subject for almost a decade
(1959--1969) -- many papers like \cite{araki}, \cite{browder} or
\cite{gugenheim} used not the final but the previous mimeographed
version: {\sc J. Milnor and J. C. Moore}, {\it On the structure of
Hopf algebras},\/ Princeton Univ., Princeton, N.J., 1959
(mimeographed)\footnote{One of the earlier pioneers told us the
  following peculiar anecdote: {\it At one point I got a hold of a
    draft copy [of Milnor--Moore paper]. A
bit later I came across an earlier draft and found it more readable,
more readable by me that is. The hunt was on for earlier drafts
because the earlier the draft, the more readable.}}
.

\section{Algebraic groups}

\subsection{Jean Dieudonn\'e}

Jean Dieudonn\'e wrote a series of papers dedicated to formal
Lie groups and hyperalgebras. In the first paper of the series
\cite{dieud54}, he starts by mentioning that:

\medbreak \noindent {\em Les r\'ecents travaux sur les groupes de
Lie alg\'ebriques
  [\ldots] montrent clairement que, lorsque le corps de base a une
  charact\'eristique $\neq 0$, le m\'ecanisme de la th\'eorie
  classique des groupes de Lie ne s'applique plus \ldots}
\medbreak

He is referring to the problems that appear when trying to extend
the dictionary `Lie groups--Lie algebras' to positive
characteristic. Dieudonn\'e was looking for the right object to
spell out an analogous dictionary. Then, in \cite[p. 87]{dieud54},
he describes the goal of the paper:

\medbreak \noindent {\em L'objet de ce travail est de d\'efinir,
pour un groupe de Lie sur
  un corps de caract\'eristique $p > 0$, une ``hyperalg\`ebre''
  associative qui correspond \`a l'alg\`ebre enveloppante de
  l'alg\`ebre de Lie dans le cas classique, mais poss\`ede ici une
  structure beaucoup plus compliqu\'ee \ldots}
\medbreak

In this paper, the author considers a formal Lie group $G$ in the
sense of \cite{bochner}, to which he attaches an associative
algebra $\mathfrak G$, called the hyperalgebra. The hyperalgebra
reflects well the structural properties of $G$.

Dieudonn\'e does not consider an abstract object,
rather some concrete algebra attached to $G$ -- essentially the
algebra of distributions with support at the identity. In
characteristic 0, this is just the universal enveloping algebra of
the Lie algebra of $G$.

The second paper of the series is devoted
to the question of when a hyperalgebra arises from a formal Lie
group, and to the study of abelian formal Lie groups. In particular
it includes the
classification of the hyperalgebras of the 1-dimensional groups.
The third and the fourth parts are also about abelian formal Lie
groups. For the theme of the present notes, the main novelty
appears in the fifth paper of this series \cite{dieud56}. In the
previous papers, the hyperalgebras in sight were mostly
commutative (because the formal groups under consideration were
commutative). Now the author wishes to pass to the general case.

\medbreak {\em Commen\c{c}ons par rappeler rapidement les
relations entre les notions d'hy\-peralg\`ebre et de groupe formel
\dots comme me l'a fait observer \textsc{P. Cartier}, la mani\`ere
la plus suggestive dont on puisse pr\'esenter cette th\`eorie
consiste \`a faire usage de la dualit\'e en alg\`ebre lineaire.}
\medbreak

He now gives a formal definition of hyperalgebra
that includes a comultiplication map; dualizing this map, he gets
the multiplication of a formal group.

\subsection{Pierre Cartier}
We discuss in this section three papers by Cartier labeled by the date
of publication.

\medbreak \textsf{1955-56.} A formal abstract definition of
\emph{hyperalgebra} is given by Pierre Cartier in \cite[Expos\'e
no. 2]{cartier-sophus}. He assumes that $K$ is a commutative ring
with unit.

\medbreak \noindent\emph{On appelle hyperalg\`ebre sur $K$ une
$K$-alg\`ebre $U$ associative, unitaire, augment\'ee, filtr\'ee,
munie d'une application diagonale v\'erifiant les axiomes
suivants.}

\begin{description}
    \item[(HA$_1$)] \emph{L'augmentation $\varepsilon: U\to K\cdot 1$ est un
    $K$-homomorphisme unitaire; on a donc $U = U^+ \oplus K\cdot 1$,
    $U^+$ designant le noyau de $\varepsilon$.}

\medbreak    \item[(HA$_2$)] \emph{L'application diagonale
$\Delta:U \to U\otimes_K U$ est un $K$-homo\-morphisme unitaire de
$U$ dans $U\otimes U$
    poss\'edant les propri\'et\'es suivantes:}
    \begin{enumerate}
        \item[a)] \emph{$I$ d\'esignant l'automorphisme identique de
        $U$, on a:}
\begin{equation}\label{cartier-1} (I\otimes \Delta)\circ \Delta = (\Delta\otimes I)\circ
\Delta.
\end{equation}
\emph{L'egalit\'e des 2 applications de $U$ dans $U\otimes
U\otimes U$ exprim\'ee par \eqref{cartier-1} s'\'enoncera en
disant que $\Delta$ est associatif. }

        \item[b)] \emph{$S$ d\'esignant l'involution $a\otimes b\mapsto b\otimes a$ de
        $U\otimes U$, on a:}
\begin{equation}\label{cartier-2} S\circ \Delta = \Delta.
\end{equation}
\emph{Autrement dit, $\Delta$ est commutative.}

        \item[c)] \emph{$P: U\otimes U \to U$ d\'esignant
        l'application $a\otimes b\mapsto ab$ d\'eduite du produit, on a:}
\begin{equation}\label{cartier-3} P\circ (I\otimes \varepsilon)\circ \Delta
= P\circ (\varepsilon\otimes I)\circ \Delta = I.
\end{equation}
\emph{Autrement dit, $\Delta$ est compatible avec l'augmentation.}
    \end{enumerate}

\medbreak    \item[(HA$_3$)] \emph{Il existe une filtration
d'alg\`ebre croissante de $U$ par des sous-modules $U_{(n)}$, $n=
0,1,\dots$, qui v\'erifie:}
\begin{align}\label{cartier-4}
U_{(0)} &= K\cdot 1
\\\label{cartier-5} \Delta(U_{(n)}) &\subset \sum_{0\le p\le n}U_{(p)} \otimes U_{(n-p)}, \qquad
\text{pour tout }n\ge 0.
\end{align}

\end{description}

Axioms (HA$_1$) and (HA$_2$), define a cocommutative bialgebra
--there is no explicit reference to the antipode. However, the
hypothesis (HA$_3$) implies the existence of the antipode. Indeed,
\eqref{cartier-5} says that $\{U_{(n)}: \, n= 0,1, \cdots \}$ is a
coalgebra filtration, hence $U_{(0)}$ is the coradical of $U$
\cite[5.3.4]{Mo}; it follows then from \eqref{cartier-4} and using
a Lemma of Takeuchi -- see \cite[5.2.10 and 5.2.11]{Mo}--  that
$U$ has a bijective antipode. \label{antipode}

\medbreak The dual of a hyperalgebra is recognized as carrying the
dual operations, so that it is -- in present terminology -- a
topological commutative bialgebra but in the paper no formal
definition is given.

Three examples of hyperalgebras are discussed in detail: the
universal enveloping algebra of a Lie algebra, the restricted
enveloping algebra of a $p$-Lie algebra and the algebra of divided
powers.

\medbreak \textsf{1956.}
 The algebra of polynomial functions on a linear
algebraic group is a commutative Hopf algebra. This fact would
seem to be the starting point of the theory of Hopf algebras to
someone with a familiarity, at least superficial, with algebraic
geometry. This is indeed recognized, albeit implicitly, in the
Comptes Rendues note \cite{cartier-cras}. Here $\Ku$ is a field.
The ultimate goal of this note is a Tannaka-type theorem.

\medbreak \noindent\emph{Soit $G$ un groupe alg\'ebrique de
matrices $g= \Vert g_{ij}\Vert$, $\mathfrak o(G)$ l'ensemble des
fonctions sur $G$ de la forme $f(g) = P(g_{ij}, \det g^{-1})$,
o\`u $P$ est un polynome. On v\'erifie facilement les
propri\'et\'es suivants:}

\begin{enumerate}
    \item[G1.] \emph{L'alg\'ebre $\mathfrak o(G)$ est une alg\'ebre de
    type fini et contient la fonction constante \'egale a 1.}

    \item[G2.] \emph{Si $g\neq g'$, il existe $f\in \mathfrak o(G)$ avec
$f(g)\neq f(g')$; un homomorphisme $\chi$ non nul de $\mathfrak
o(G)$ dans $\Ku$ est de la forme $\chi(f) = f(g_0)$ $(g_0\in G)$.}

    \item[G3.] \emph{Si $f\in \mathfrak o(G)$ et $f'(g_1, g_2) =
    f(g_1^{-1}g_2)$, on a $f'\in \mathfrak o(G) \otimes \mathfrak o(G)$.}
\end{enumerate}

The convolution product in the space of linear forms on $\mathfrak
o(G)$ is defined, and the hyperalgebra $\mathbf U (G)$ is
described as the space of those linear forms vanishing on some
power of the augmentation ideal of $\mathfrak o(G)$. The existence
of the comultiplication of $\mathbf U(G)$ is made explicit, in
terms of the product of $\mathfrak o (G)$.

The subject was continued by Gerhard Hochschild and Dan
Mostow in a series of papers that will be considered later.

\medbreak \textsf{1962.} A few years after the papers already
discussed, Cartier presented a much more ample set of ideas in
\cite{cartier-bru}. Although he recognizes that

\medbreak \noindent\emph{Le pr\'esent article n'est qu'un
r\'esum\'e \dots Nous publierons par ailleurs un expos\'e
detaill\'e de la th\'eorie.}

\medbreak \noindent(exposition that never saw the light), the
paper contains many important concepts. Let us first select some
parts of the Introduction, to highlight the main motivations of
his work. He begins by saying that the theory of algebraic groups
in positive characteristic has complications arising from
phenomena of inseparability, and mentions explicitly the problem
of isogeny, to which Barsotti and he have devoted some work.
Then he recognizes that the language of schemes -- freshly
introduced by Grothendieck at that time -- is the natural setting
to treat these questions and explains that he will proceed within
it with some restrictions. Then he mentions the work of
Dieudonn\'e on formal groups and declares:

\medbreak \noindent\emph{Mais jusq'\`a pr\'esent, la th\'eorie des
groupes alg\'ebriques et celle des groupes formels n'avaient que
des rapports superficiels. C'est le but du pr\'esent expos\'e d'en
faire la synth\`ese. Pour cela, nous devrons \'elargir un peu la
notion de groupe formel \dots ceci fait,  la th\'eorie des groupes
formels est \'equivalente \`a celle des hyperalg\`ebres.}
\medbreak

He immediately points out in a footnote:

\medbreak \noindent\emph{En topologie alg\'ebrique, on utilise
sous le nom d'Alg\`ebre de Hopf une notion tr\`es proche de celle
de hyperalg\`ebre.}

\medbreak He refers next to the duality between commutative algebraic groups 
and commutative formal groups and observes in a footnote:

\medbreak \noindent\emph{Gabriel utilise le langage des
hyperalg\`ebres, et non celui des groupes formels, ce qui oblige
\`a de nombreuses contorsions.}

\medbreak This last remark, together with the allusions to the
theory of schemes, explains the choice that Cartier makes for the
exposition. Namely, he rather prefers a systematic use of a
functorial presentation, than \emph{le langage des
hyperalg\`ebres}; Hopf algebras appear only marginally in the
paper.

\medbreak Let us now describe the main points of the paper. He
begins in Section 2 by identifying the group $GL(n, A)$, $A$ a
commutative ring, with the points of $A^r$, $r= 1 + n^2$,
annihilated by the polynomial $D = 1 - X_0 \det X_{ij}$. He then
observes that, given a field $\Omega$, an algebraic subgroup of
$GL(n, \Omega)$ is determined by an ideal $I$ of $\Omega[X_0,
X_{ij}]= \Omega[X]$ satisfying the following conditions.

\begin{enumerate}
\renewcommand{\theenumi}{\alph{enumi}}   \renewcommand{\labelenumi}{(\theenumi)}

\item $D$ belongs to $I$.

\item If $Y_0$, $Y_{ij}$ is another set of $r$ variables, abbreviated $Y$, then set
$Z_0 = X_0Y_0$, $Z_{ik} = \sum_j X_{ij}Y_{jk}$. For any $P\in I$,
there exist $L_{\alpha}, L'_{\beta}\in \Omega[X]$, $P_{\alpha},
P'_{\beta}\in I$ such that
\begin{equation}\label{cartier62-dos}
P(Z) = \sum_{\alpha}L_{\alpha}(X) P_{\alpha}(Y) + \sum_{\beta}
P'_{\beta}(X) L'_{\beta}(Y).
\end{equation}

\item Let $X'_0$ be the determinant of the matrix $(X_{kl})$ and
$X'_{ij}$ the product of $X_0$ by the cofactor of $X_{ij}$ in the
determinant $X'_0$. Then $P(X')$ belongs to $I$ whenever $P\in I$.

\item The ideal $I$ is an intersection of prime ideals.
\end{enumerate}

In other words, condition (b) says that $I$ is a coideal of
$\Omega[X]$ and (c) that it is stable under the antipode. Hence,
the quotient algebra $\Omega[X]/I$ of polynomial functions on the
algebraic subgroup has a natural structure of Hopf algebra
(without nilpotents in accordance to (d)).

\medbreak He then proceeds to define  algebraic groups as
representable functors from the category of commutative algebras
to the category of groups and devotes some time to the study of
their properties.

\medbreak Next, he defines in Section 10 the notion of formal
group also as a functor from the category of commutative algebras
to the category of groups, but instead of representable he asks
for another set of axioms that we will not recall here. These
axioms allow him to define the hyperalgebra $U(G)$ associated to a
formal group $G$; here the coproduct is explicitly mentioned,
otherwise he refers to \cite{cartier-sophus}. This assignment
gives an equivalence between the category of formal groups and the
category of hyperalgebras.

\medbreak He next states a structure result for formal groups
(Theorem 2): any formal group over a perfect field can be
decomposed in a unique way as a semidirect product of a
``separable'' group and an ``infinitesimal'' group. This
translates to a structure theorem for hyperalgebras (that is,
cocommutative Hopf algebras) that is not explicitly spelled out;
the explicit formulation of this structure theorem was discovered
by Kostant, see page \pageref{kostant-structure} below. However,
Theorem 3 says that the hyperalgebra of an infinitesimal group
over a field of characteristic 0 is a universal enveloping
algebra. Contrary to Theorem 5.18 in \cite{milnormoore} quoted
above, it is \emph{not} assumed that the algebra is primitively
generated.

\subsection{Gerhard Hochschild and  Dan Mostow} These
authors start the Introduction
of the first of the series of three
papers that appeared in the late 50's \cite{representative1},
\cite{representative2} and \cite{representative3} as follows:

\medbreak \noindent{\it We are concerned with a number of
inter--related general
  questions concerning the representations of a Lie group $G$ and the
  algebra $R(G)$ of the (complex valued) representative functions on
  $G$ \ldots
Instead of dealing with the characters of $R(G)$, we deal with its proper
  automorphisms [A = group of all algebra automorphisms commuting with right
  translations] \ldots The left translations on $R(G)$ [extend] to an
  operation of $G^+$ [\, the universal complexification of G\,]
on $R(G)$, so that one obtains a continuous monomorphism
  of $G^+$  into $A$ \ldots Chevalley's complex formulation of
  Tannaka's theorem becomes the statement that, if $G$ is compact,
  this monomorphism sends $G^+$ onto $A$. [ If $R(G)$ is finitely
  generated] \, we may identify $A$ and hence $G^+$ with an algebraic
  complex linear group. This links up the Lie group situation with the
analogous situation for algebraic groups, which is dealt with by
P. Cartier [\,in \cite{cartier-cras}\,].}

\medbreak

In the second and third of the series of papers, the same
platform is considered but under
hypotheses for $R(G)$ more general than the finite generation.

For example, in the first paper it is proved that if $R(G)$ is
finitely generated, then the image of the natural morphism from
$G^+$ into $A$ is $A$ itself. This is not true in general and in
\cite{representative2} this image is explicitly calculated.

\medbreak \noindent{\em Our result is that a proper automorphism
$\alpha$ belongs to
  the natural image of $G^+$ if and only if $\alpha(\operatorname{exp}(h))=
  \operatorname{exp}(\alpha(h))$ for every $h \in
  \operatorname{Hom}(G,C)$}.
\medbreak

In the third paper \cite{representative3} the authors deal with
complex groups. In their own words:

\medbreak \noindent{\em The theory developed \,[so far\,] was
adapted to real
  Lie groups, and used complex Lie groups as auxiliaries. It is evident
  from the nature of the proofs and the main results that there is an
  underlying analogous theory for complex Lie groups, which is
  actually simpler than the theory for real Lie groups. Nevertheless,
  the systematic development of the complex analytic case brings up a
  number of new questions which are of independent interest.}

\medbreak

Constantly, along this series of papers, the authors work --
without explicitly spelling out the words -- with the natural Hopf
algebra structure on $R(G)$. The algebra structure appears in an
explicit form but the coalgebra structure appears under the guise
of the formul\ae\/ that are obtained -- once we know that for $f
\in R(G)$ the set $G\cdot f$ generates a finite dimensional
subspace -- when we express $x \cdot f = \sum \alpha_i(x) g_i$ in
terms of a basis  $\{g_1,\cdots , g_n\}$ of this subspace.  In
terms of the comultiplication of $R(G)$ the above formula just
means that $\Delta(f)=\sum g_i \otimes \alpha_i$.

\medbreak

It might be important to stress at this point that this series of
papers deals globally with a generalization of the
Tannaka theorem and in that sense, follow the same track initiated by
Cartier in \cite{cartier-cras} and that was considered before. Here again
we are observing a fact that is very relevant even today:
the deep interaction existing between Tannaka reconstruction
viewpoint and Hopf theory.

\medbreak Moreover, the work of Hochschild from the late 50's
until the late 70's contained a systematic approach to the theory
of Lie and algebraic groups and their Lie algebras, with the
viewpoint of Hopf theory in mind. As time passed, the Hopf
structure started to appear explicitly and its use became more and
more systematic. In particular in the book {\em The structure of Lie
groups} (1965), \cite{hoch-65}, the Hopf algebra structure plays
in many occasions a conspicuous role. For example, on page 26 of
Chapter 2, $\S$ 3 of the mentioned reference, the author considers
the $\R$--algebra $R(G)$ of real valued representative functions
on a compact group $G$ and proves that it admits a natural
structure of ``Hopf algebra''.

\medbreak \noindent{\it If we abstract the structure we have just
described from
  the group $G$, retaining only the $\R$--module ${\mathcal R}(G)$,
  the algebra structure $\mu$ on ${\mathcal R}(G)$ with the algebra
  homomorphism $u:\R \rightarrow {\mathcal R}(G)$ as a unit and the
  coalgebra structure $\gamma$ with the algebra homomorphism
  $c: {\mathcal R}(G) \rightarrow \R$ as a counit, and if we assume
  that the above identities are satisfied by these maps and that
  $\gamma$ is an algebra homomorphism, we have what is called a Hopf
  algebra $({\mathcal R}(G), \mu, u, \gamma, c)$.}
\medbreak

The reader should notice that the antipode was not present in the
above definition of Hopf algebra. But the author goes on and says
on page 27 op.cit.:

\medbreak \noindent{\it In our special case, there is one further
structural
  item. Let $\eta$ denote the algebra endomorphism of ${\mathcal
    R}(G)$ that is defined from the inverse of $G$, i.e.,
  $\eta(f)=f'$, where $f'(x)=f(x^{-1})$. Then one sees immediately
  that $\mu \circ (\eta \otimes i) \circ \gamma = u \circ c$. Let us
  agree to call an algebra endomorphism $\eta$ with this last property
  a symmetry of our Hopf algebra.}
\medbreak

It is worth noticing that in that same period, the Hopf algebra
language in algebraic group theory was becoming more and more
standard --see for example the work by Demazure, Gabriel and
Grothendieck, \cite{dga},\cite{dg1}.

\medbreak We finish by saying some words about the evolution of
the concept of integral in the setting of Hopf algebras. Since
Hurwitz and Weyl, invariant integration was one of the main tools
in representation theory and was formalized in the context of Hopf
algebra theory in the early sixties.

The fist appearance of the general concept of integral in the
framework of Hopf theory seems to be in \cite[page 29]{hoch-65}
under the name of a ``gauge''.

\medbreak \noindent{\it Now suppose that $G$ is compact. Then a
Haar integral
  for $G$ gives us an $\R$ module homomorphism $J: {\mathcal R}(G)
  \rightarrow \R$. The invariance property $J(y\cdot f)=J(f)$ for all
  $y \in G$ and all $f \in {\mathcal R}(G)$, is easily seen to be
  equivalent to the relation $(J \otimes i)\circ \gamma = u \circ J$
  $\ldots$ If [$\,A\,$] is any Hopf algebra let us call a map $J:A
  \rightarrow \R$ with [\,this\,] property a gauge of the Hopf algebra.}
\medbreak

The relationship between the existence of a ``gauge'' in the Hopf algebra
of functions over a group  and the semisimplicity of the representations
of the latter, appears clearly
in the above mentioned book in the guise of the following {\em Tannaka
duality theorem} (page 30):

\medbreak \noindent{\it For a compact group $G$, let $\mathcal
H(G)$ denote the
  Hopf algebra attached to $G$ in the canonical fashion. Then
  $\mathcal H(G)$ is a reduced Hopf algebra with symmetry and
  gauge. For a reduced Hopf algebra having a symmetry and a gauge, let
$\mathcal G(H)$ be the topological group of homomorphisms attached
to $H$ in the canonical fashion. Then $\mathcal G(H)$ is compact.
The canonical maps $G \rightarrow \mathcal {\mathcal G}(\mathcal
H(G))$ and $H \rightarrow \mathcal H(\mathcal G(H))$ are
isomorphisms.} \medbreak

The same concept of gauge is also considered later in the context
of Lie theory -- see \cite{hoch-70} -- where in Theorem 6.2 it is
proved that:

\medbreak \noindent{\it Let $L$ be a finite dimensional Lie
algebra over the field F of characteristic 0. Then $L$ is
semisimple if and only if there is a linear functional $J$ on
${\bf{H}}(L)$ such that $J \circ u$ is the identity map on $F$,
and $u \circ J= (J \otimes i)\circ \gamma$ where $u$ is the unit
of ${\bf{H}}(L)$ and $\gamma$ is the comultiplication.} \medbreak

Here ${\bf{H}}(L)$ is the algebra of linear functions on the enveloping algebra
$U(L)$ that vanish on a two-sided ideal of finite codimension.
Later, these ideas were generalized by M. Sweedler in
\cite{Swintegrals} to the context of arbitrary Hopf algebras:

\medbreak \noindent{\it For a Hopf algebra that is the coordinate
ring of a
  compact Lie group, there is a unique one--dimensional left invariant
  ideal in the linear dual, this is the space spanned by a left
Haar integral. Hochschild has observed that for a Hopf algebra $H$
which is the coordinate ring of an affine algebraic group the
existence of a left invariant ideal in the linear dual which
complements the linear functionals orthogonal to the unit of $H$, is
equivalent to the group being completely reducible.

In this paper \ldots [\,we\,] generalize Hochschild's results.}

\medbreak Integrals in finite-dimensional Hopf algebras were
previously considered in \cite{LS}, cf. also \cite{kp66}. The
notion of integral plays a central role in the definition of ring
group by Georg I. Kac discussed below.

\subsection{Bertram Kostant}\label{kostant} Bertram Kostant,
besides being the thesis advisor to Moss Sweedler and the
originator of some of the current notations and nomenclature in
Hopf theory -- for example of the term \emph{group-like element}
-- produced one of the first genuine applications of Hopf algebras
with his contribution to the understanding of finite groups of Lie
type, in his remarkable paper \cite{Ko-Z}.

\medbreak Here is the limpid definition of a Hopf algebra given in
the Preliminaries of \cite{Ko-Z}. Kostant begins by introducing
the convolution product $\ast$ in $\Hom (A,R)$, where $A$ is a
coalgebra (with comultiplication $d$ and counit $\varepsilon$)
over a commutative ring $C$, and $R$ is an algebra (with
multiplication $m$ and unit $\rho$) still over $C$. Then he says:

\medbreak
\noindent
\emph{Now assume that $A$ is a Hopf algebra
($A$ is an algebra and coalgebra such that $d$ and $\varepsilon$ are
homomorphisms and $\varepsilon\rho$ is the identity on $C$).
\newline
By an antipode on $A$ we mean an element (necessarily unique if it
exists) $s\in \Hom_C(A,A)$ such that $I\ast s = s\ast I =
\varepsilon$ where $I$ is the identity on $A$ and $\ast$ is as
above with $A$ taken for $R$. From now on Hopf algebra means Hopf
algebra with antipode\footnote{The referee asked us if this was
the first time the term \emph{antipode} was used. We do not know;
it was named a \emph{symmetry} in \cite{hoch-65} and an
\emph{antipodisme} in \cite[p. 511]{dg1}. Mariano Su\'arez \'Alvarez pointed out to us
that the word \emph{antipodism} is used in \cite[Chap. XI, Sect. 8]{CE}}.}

\medbreak He continues by observing that a Hopf algebra gives rise
to a representable functor from the category of algebras over $C$
to the category of groups.

Still in the Preliminaries, given a Hopf algebra $B$ over $\mathbb
Z$, he explains how to construct  from any `admissible' family of
ideals of finite type a  Hopf algebra dual to $B$.

The bulk of the paper contains the definition of a `divided
powers' version over $\Z$ of the universal enveloping algebra of a
finite-dimensional semisimple Lie algebra, i.e., the so called
``Kostant $\Z$--form''. In accordance with a private communication
by the author:

\medbreak {\em The ultimate purpose of this [...] was to exhibit
the Chevalley (Tohoku) finite groups as the group-like elements in
the $\Z$ dual of $U_{\mathbb Z}(\mathfrak g)$. However I did not
succeed in doing this. However recently George Lusztig, in a tour
de force paper, did succeed in establishing what I hoped would be
the case. Lusztig's paper can be found in ArXiv:0709.1286v1
[math.RT] 9 Sept. 2007\footnote{{\em Study of a $ \mathbf Z$-form of the
  coordinate ring of a reductive group}, to appear in J. Amer. Math. Soc.}.}

\medbreak Kostant also proved two theorems on the structure of
cocommutative Hopf algebras over an algebraically closed field
$k$.\label{kostant-structure} First, any cocommutative Hopf
algebra $H$ is the smash product of an irreducible Hopf algebra
$U$ with a group algebra $k G$. Here $G$ turns out to be the group
$G(H)$ of group-likes in $H$, and `irreducible' means that $G(U)$
is trivial. Second, $U$ is the universal enveloping algebra of the
Lie algebra of primitive elements in $H$, provided that the
characteristic of $k$ is 0. These results were never published by
Kostant\footnote{Susan Montgomery kindly pointed out to us that,
according with some correspondence with Sweedler or Larson,
Kostant proved this result in the early 1960's.}, and appeared
first in the book \cite{Sw}, see Theorem 8.1.5 and Section 13.1.

\medbreak  Later, Kostant obtained the analogous result for
cocommutative Hopf super algebras \cite[Th. 3.3]{Ko-super}. In
this paper he introduced the notion of super manifolds and super
Lie groups with the goal to do in the super case, what Chevalley
did for ordinary Lie groups. As the author told us in a private
communication: {\em
  the idea is to get the super Lie group from its Lie algebra using
  the Hopf algebra double dual}.

\subsection{Georg Isaakovich Kac}

Inspired by work of Stinespring  on a duality theory for
unimodular groups (extending previous work by Tannaka and Krein),
Georg I. Kac introduced in 1961, see \cite{Kac}, the notion of
\emph{ring group}. We quote from the introduction:

\medbreak \noindent\emph{The basic idea of the generalization of
the concept of group consists in the following. Let $\mathfrak G$
be a group. We denote by $\mathfrak M$ the commutative ring (with
respect to the operation of multiplication) of bounded functions
on $\mathfrak G$. The mapping $f(y) \to f(xy)$ is an isomorphism
of $\mathfrak M$ into $\mathfrak M\otimes \mathfrak M$. Thus
$\mathfrak G$ can be considered as a commutative ring $\mathfrak
M$ with a given isomorphism of $\mathfrak M$ into $\mathfrak
M\otimes \mathfrak M$. We obtain the generalization under study,
if we reject in this definition the requirement that $\mathfrak M$
is commutative.}

\medbreak A ring group is essentially like a Hopf algebra in the
setting of von Neumann algebras. Namely, it is defined as a
collection $(\mathfrak M,\Phi,{}^+, m)$, where $\mathfrak M$ is a
von Neumann algebra; $\Phi$ is a monomorphism of $\mathfrak M$
into $\mathfrak M\otimes \mathfrak M$\footnote{Kac says
``isomorphism\ldots into'' which meant at the
  time ``monomophism''; for isomorphisms in our present meaning, it
  was used ``isomorphism\ldots onto''.},
that is coassociative; $A\mapsto A^+$ is a sesquilinear
involution of $\mathfrak M$ that reflects somehow the inversion in
the group; $m$ is a measure on $\mathfrak M$ (a `gauge'),
satisfying certain axioms, similar to what is called today an
integral.

\medbreak A ring group is representable as the ring of bounded
functions on a unimodular group $G$ if and only if the ring
$\mathfrak M$ is commutative. Every ring group $\mathfrak M$
corresponds uniquely to a so-called dual Hilbert ring; this
formalism enables one to define the dual $\mathfrak M\,\hat{}$ of
$\mathfrak M$. There is a duality theorem stating that the double
dual $\mathfrak M\,\hat{}\,\,\hat{}$ \,is isomorphic to $\mathfrak
M$.

\medbreak Georg I. Kac developed the theory of finite ring groups
(ring groups of finite dimension, which are always semisimple algebras)
in several papers. For instance,
in \cite{kac62}, he extends certain classical finite-group
theorems to finite-ring groups, e. g. the analogue of Lagrange's
theorem: The order of a subgroup of a ring group divides the order
of the ring group --this was generalized in 1989 to finite
dimensional Hopf algebras over arbitrary fields by Nichols and Z\"oller --.
In \cite{kp66, kac68}, the construction of new
finite ring groups by extensions is given; this was rediscovered
by Takeuchi in 1981 and then again by Majid in 1988. In
\cite{kp66}, the existence of an integral in a finite ring group
is derived, assuming the existence of an antipode; this was later
established for finite-dimensional Hopf algebras by Larson and
Sweedler \cite{LS}. In \cite{kac72}, among other results, he
showed that a finite ring group of prime dimension is necessarily
a group algebra. This fact, for general Hopf algebras over $\mathbb C$
of dimension
$p$, was conjectured by Kaplansky in 1975
and proved by Y. Zhu in 1993.

\section{Final considerations}


With the publication of Sweedler's book \cite{Sw} -- September 1,
1969 -- the subject started to shape up as an independent part of
abstract algebra and was able to walk without the help of its two
honorable parents. Of the latter period we only mention that, with
the appearance in 1987 of the paper by V. Drinfel'd
\cite{drinfeld} and the subsequent work by him and many other
mathematicians, the area experienced a very radical change in
terms of methods, examples and interaction with other parts of
mathematics. In this latter period, the developments obtained in
understanding the structure of Hopf algebras and their
representations have been outstanding, and they have been entwined
with the development of different areas of mathematics, most
remarkably with: knot theory and topology, conformal field theory,
ring theory, category theory, combinatorics, etc.

\end{document}